\documentclass[12pt,dvips]{amsart}
\usepackage{amsfonts, amssymb, latexsym, epsfig}

\newtheorem{Theorem}{Theorem}[section]

\newtheorem{Main Conjecture}[Theorem]{Main Conjecture}

\theoremstyle{remark}
\newtheorem{Example}[Theorem]{Example}

\newcommand{\excise}[1]{}

\begin{document}
\pagestyle{plain}
\title{Noncoherent initial ideals in exterior algebras}
\author{Dominic Searles}
\address{Dept. of Mathematics\\
University of Illinois at Urbana-Champaign\\
Urbana, IL 61801}
\email{searles2@illinois.edu}

\author{Arkadii Slinko}
\address{Dept. of Mathematics\\
University of Auckland\\
Auckland, New Zealand}
\email{a.slinko@auckland.ac.nz}

\date{\today}

\begin{abstract}
We construct a noncoherent initial ideal of an ideal in the exterior algebra of order $6$, answering a question of D. Maclagan (2000). We also give a method for constructing noncoherent initial ideals in exterior algebras using certain noncoherent term orders.
\end{abstract}

\maketitle

\section{Introduction}

Let $\mathbb{F}$ be an algebraically closed field of characteristic $\neq 2$. The {\bf exterior algebra} $\bigwedge$ of order $n$ over $\mathbb{F}$ consists of polynomials with coefficients from $\mathbb{F}$ in noncommuting indeterminates $x_1,\ldots, x_n$ subject to the relation $x_ix_j = -x_jx_i$ for $1\le i,j\le n$. This relation implies $x_i^2=0$ and that any monomial can be reordered up to a sign change into a canonical form $x_{i_1}\ldots x_{i_k}$ where $i_1 < \ldots < i_k$. Throughout, let $\mathbf{x}^{\mathbf{a}}$ denote $x_1^{a^{(1)}}\ldots x_n^{a^{(n)}}$, where $\mathbf{a}=(a^{(1)},\ldots,a^{(n)})\in \mathbb{N}^n$. Then $M_{\bigwedge} = \{\mathbf{x}^{\mathbf{a}}:\mathbf{a}\in\{0,1\}^n\}$ is the set of monomials of $\bigwedge$.

A {\bf term order} $\prec$ on $M_{\bigwedge}$ is a total order on $M_{\bigwedge}$ which satisfies:

\begin{itemize}
\item[(i)]  $1= \mathbf{x}^{\mathbf{0}} \prec \mathbf{x}^{\mathbf{a}}$ for all $\mathbf{x}^{\mathbf{a}} \neq 1$ in $M_{\bigwedge}$.
\item[(ii)] If $\mathbf{x}^{\mathbf{a}} \prec \mathbf{x}^{\mathbf{b}}$ then $\mathbf{x}^{\mathbf{a}+\mathbf{c}}\prec \mathbf{x}^{\mathbf{b}+\mathbf{c}}$ whenever $\mathbf{x}^{\mathbf{a}+\mathbf{c}}$ and $\mathbf{x}^{\mathbf{b}+\mathbf{c}}$ are both in $M_{\bigwedge}$.
\end{itemize}

A term order $\prec$ on $M_{\bigwedge}$ is {\bf coherent} if there exists a weight vector $w\in \mathbb{R}^n$ such that $w\cdot \mathbf{a} < w\cdot \mathbf{b}$ whenever $\mathbf{x}^{\mathbf{a}}\prec \mathbf{x}^{\mathbf{b}}$, and {\bf noncoherent} otherwise. Equivalently, a term order on $M_{\bigwedge}$ is coherent if it can be extended to a term order on the monomials of the usual (commutative) polynomial algebra $\mathbb{F}[x_1,\ldots,x_n]$. When $n\ge 5$, there exist noncoherent term orders on $M_{\bigwedge}$. In the interpretation of term orders on $M_{\bigwedge}$ as comparative probability orders on subsets of an $n$-element set, this fact has long been known (\cite{Kraft.Pratt.Seidenberg}). 

Let $(\mathbf{x}^{\mathbf{a}_1}, \ldots , \mathbf{x}^{\mathbf{a}_N}) =_0 (\mathbf{x}^{\mathbf{b}_1}, \ldots , \mathbf{x}^{\mathbf{b}_N})$ mean that $\mathbf{x}^{\mathbf{a}_i}, \mathbf{x}^{\mathbf{b}_i} \in M_{\bigwedge}$ for all $1\le i\le N$ and $\sum_{1\le i\le N} \mathbf{a}_i = \sum_{1\le i\le N} \mathbf{b}_i$. The following condition, when added to the defining conditions (i), (ii) of a term order on $M_{\bigwedge}$, results in a set of conditions which are both necessary and sufficient for coherency: 
\begin{itemize}
\item[(iii)] For all $N \ge 2$ and all $\mathbf{x}^{\mathbf{a}_i}, \mathbf{x}^{\mathbf{b}_i} \in M_{\bigwedge}$, if $(\mathbf{x}^{\mathbf{a}_1}, \ldots , \mathbf{x}^{\mathbf{a}_N}) =_0 (\mathbf{x}^{\mathbf{b}_1}, \ldots , \mathbf{x}^{\mathbf{b}_N})$ and $\mathbf{x}^{\mathbf{a}_i}\prec \mathbf{x}^{\mathbf{b}_i}$ for all $i < N$, then it is not the case that $\mathbf{x}^{\mathbf{a}_N}\prec \mathbf{x}^{\mathbf{b}_N}$.
\end{itemize}
The equivalent formulation of this condition in the setting of comparative probability orders can be found in, e.g., \cite{Fishburn96}. Violation of (iii) for a certain $N$ is called a failure of the $N^{th}$ {\bf cancellation condition}, denoted $C_N$, and implies the order is noncoherent.   

Let $\prec$ be a term order and let $f = \sum_i c_i \mathbf{x}^{\mathbf{a}_i} \in \bigwedge$, where $0\neq c_i\in \mathbb{F}$ and $\mathbf{x}^{\mathbf{a}_i}\in M_{\bigwedge}$. Then the {\bf initial monomial} $in_\prec(f)$ of $f$ is $\max_i \mathbf{x}^{\mathbf{a}_i}$ (where the maximum is taken with respect to $\prec$), and the {\bf lead term} $LT_\prec(f)$ of $f$ is $\max_i c_i \mathbf{x}^{\mathbf{a}_i}$. Let $I\subseteq \bigwedge$ be a left ideal. Then the {\bf initial ideal} $in_\prec(I)$ of $I$ with respect to $\prec$ is the monomial ideal left-generated by $\{in_\prec(f):f\in I\}$. An initial ideal $in_\prec(I)$ with respect to some noncoherent term order $\prec$ is a {\bf noncoherent initial ideal} if $in_\prec(I)\neq in_{\prec_c}(I)$ for any coherent term order $\prec_c$. We will work only with homogeneous ideals $I\subset \bigwedge$. These are in fact two-sided, so we may drop the word ``left'' from the discussion (see, e.g., \cite[Section 7]{Stokes}). 

In \cite{Maclaganthesis}, D.~Maclagan posed the following question:

\begin{quote}
Does there exist a \emph{noncoherent initial ideal} of an ideal $I$ in the exterior algebra? 
That is, is there some initial ideal of $I$ with respect to some noncoherent term order which is not equal to the initial ideal of $I$ with respect to any coherent term order?
\end{quote}

We give an affirmative answer to this question. For $\mathbf{x}^{\mathbf{a}}\in M_{\bigwedge}$, let $|\mathbf{a}|$ denote the sum of the entries of $\mathbf{a}$. Let $\mathcal{P}$ denote the lattice of elements of $\{0,1\}^n$, in which $\mathbf{a}<\mathbf{b}$ whenever $\mathbf{b}-\mathbf{a}\in\{0,1\}^n$ (this is the same as ordering the subsets of an $n$-element set by inclusion). By an {\bf antichain} in this lattice, we mean a set of pairwise incomparable elements.

\begin{Theorem}\label{theorem:main}
Suppose $\prec$ is a noncoherent term order with a $C_N$ failure: $(\mathbf{x}^{\mathbf{a}_1}, \ldots , \mathbf{x}^{\mathbf{a}_N}) =_0 (\mathbf{x}^{\mathbf{b}_1}, \ldots , \mathbf{x}^{\mathbf{b}_N})$ and $\mathbf{x}^{\mathbf{a}_i}\prec \mathbf{x}^{\mathbf{b}_i}$ for all $1\le i\le N$. If $|\mathbf{a}_i|=|\mathbf{b}_i|$ for all $1\le i \le N$, and $\{\mathbf{a}_1, \ldots , \mathbf{a}_N, \mathbf{b}_1, \ldots , \mathbf{b}_N\}$ is an antichain in $\mathcal{P}$, then the initial ideal $in_{\prec}(I)$ of the (homogeneous) ideal 
\[I = \langle \, \mathbf{x}^{\mathbf{b}_i} - \mathbf{x}^{\mathbf{a}_i}, 1 \le i \le N \, \rangle\subset \bigwedge \]
is a noncoherent initial ideal.
\end{Theorem} 

It is not obvious that there exists an order $\prec$ satisfying the hypotheses of Theorem~\ref{theorem:main}. We exhibit an example of such an order, which was found using the MAGMA computer algebra system \cite{MAGMA}.

\begin{Example}\label{ex:main}
Let $\prec$ denote the following term order on the exterior algebra of order $6$:

\noindent
$1 \prec x_1 \prec x_2 \prec x_3 \prec x_1x_2 \prec x_1x_3 \prec x_4 \prec x_5 \prec x_1x_4 \prec x_6 \prec x_2x_3 \prec x_1x_5 \prec x_1x_6 \prec x_1x_2x_3 \prec x_2x_4 \prec x_2x_5 \prec x_3x_4 \prec x_1x_2x_4 \prec x_3x_5 \prec x_2x_6 \prec \boxed{x_1x_2x_5 \prec x_1x_3x_4} \prec x_3x_6 \prec \boxed{x_1x_3x_5 \prec x_1x_2x_6} \prec x_4x_5 \prec x_1x_3x_6 \prec x_4x_6 \prec \boxed{x_2x_3x_4 \prec x_1x_4x_5} \prec x_5x_6 \prec \boxed{x_1x_4x_6 \prec x_2x_3x_5}\prec \ldots$ 

(the remaining comparisons are determined since $\mathbf{x}^{\mathbf{a}}\prec\mathbf{x}^{\mathbf{b}} \iff \mathbf{x}^{\mathbf{1-b}}\prec\mathbf{x}^{\mathbf{1-a}}$ when $\mathbf{x}^{\mathbf{a}},\mathbf{x}^{\mathbf{b}}\in M_{\bigwedge}$).

The four boxed comparisons
\[\mathbf{x}^{\mathbf{a}_1} =  x_1x_2x_5 \prec x_1x_3x_4 = \mathbf{x}^{\mathbf{b}_1}, \qquad \mathbf{x}^{\mathbf{a}_2} = x_1x_3x_5 \prec x_1x_2x_6 = \mathbf{x}^{\mathbf{b}_2},\] \[\mathbf{x}^{\mathbf{a}_3} = x_2x_3x_4 \prec x_1x_4x_5 = \mathbf{x}^{\mathbf{b}_3}, \qquad \mathbf{x}^{\mathbf{a}_4} = x_1x_4x_6 \prec x_2x_3x_5 = \mathbf{x}^{\mathbf{b}_4}\]
are a failure of $C_4$ satisfying the hypotheses of Theorem~\ref{theorem:main}. Let $I = \langle \,  \mathbf{x}^{\mathbf{b}_i}-\mathbf{x}^{\mathbf{a}_i}:1\le i\le 4 \, \rangle$. Using the Gr\"obner basis algorithm of \cite[Theorem 6.6]{Stokes}, we compute 
\[in_{\prec}(I) = \langle \, \mathbf{x}^{\mathbf{b}_1}, \,\,  \mathbf{x}^{\mathbf{b}_2}, \,\, \mathbf{x}^{\mathbf{b}_3}, \,\, \mathbf{x}^{\mathbf{b}_4}, \,\, x_1x_3x_5x_6, \,\, x_2x_3x_4x_6 \, \rangle.\]
It is actually easy to see that $in_\prec(I)$ is noncoherent. By condition (iii), at least one of $\mathbf{x}^{\mathbf{a}_1}, \mathbf{x}^{\mathbf{a}_2}, \mathbf{x}^{\mathbf{a}_3}, \mathbf{x}^{\mathbf{a}_4}$ must appear in $in_{\prec_c}(I)$ for any coherent term order $\prec_c$. It is clear that none of these monomials are in $in_\prec(I)$.
\end{Example}

\section{Proof of Theorem~\ref{theorem:main}}

For $\mathbf{x}^{\mathbf{a}},\mathbf{x}^{\mathbf{b}}\in M_{\bigwedge}$, we say $\mathbf{x}^{\mathbf{a}}$ {\bf divides} $\mathbf{x}^{\mathbf{b}}$ if $\mathbf{b-a}\in \{0,1\}^n$, and by $\frac{\mathbf{x}^{\mathbf{b}}}{\mathbf{x}^{\mathbf{a}}}$ we mean the monomial $\mathbf{x}^{\mathbf{b-a}}$. By the {\bf least common multiple} of $c_i\mathbf{x}^{\mathbf{a}}$ and $c_j\mathbf{x}^{\mathbf{b}}$, we mean $lcm(c_i,c_j)\cdot\mathbf{x}^{\mathbf{a}\cup\mathbf{b}}$, where $(\mathbf{a}\cup\mathbf{b})^{(i)}=0$ if $\mathbf{a}^{(i)}=\mathbf{b}^{(i)}=0$, and $(\mathbf{a}\cup\mathbf{b})^{(i)}=1$ otherwise.

Fix a term order $\prec$ and let $G=\{ g_1, \ldots , g_r \}\subset \bigwedge$. Suppose $\mathbf{x}^{\mathbf{c}}  in_{\prec}(g_i) = 0$ for some $g_i\in G, \mathbf{x}^{\mathbf{c}}\in M_{\bigwedge}$. Following the notation of \cite{Maclaganthesis}, define the {\bf $T$-polynomial} $T_{g_i, \mathbf{x}^{\mathbf{c}} }$ to be the polynomial 
\[T_{g_i, \mathbf{x}^{\mathbf{c}} }=\mathbf{x}^{\mathbf{c}}  g_i.\]

Let $m_{g_i,g_j}$ be the least common multiple of $LT_\prec(g_i)$ and $LT_{\prec}(g_j)$. Then the {\bf $S$-polynomial} $S_{g_i,g_j}$ is the polynomial 
\[S_{g_i,g_j}= (-1)^{d_1}\frac{m_{g_i,g_j}}{LT_{\prec}(g_i)}g_i - (-1)^{d_2}\frac{m_{g_i,g_j}}{LT_{\prec}(g_j)}g_j.\]
where $d_1=1$ if reordering $\frac{m_{g_i,g_j}}{LT_{\prec}(g_i)}in_\prec(g_i)$ into canonical form changes the sign of this monomial, $d_1=0$ otherwise, and $d_2$ is defined similarly.

Let $f\in \bigwedge$. By {\bf reducing} $f$ with respect to $G$ and $\prec$, we mean choosing some $g_i\in G$ such that $in_\prec(g_i)$ divides $in_\prec(f)$, letting $r= f-(-1)^d\frac{LT_\prec(f)}{LT_\prec(g_i)}g_i$ (where $d=1$ if reordering $\frac{in_\prec(f)}{in_\prec(g_i)}in_\prec(g_i)$ changes the sign of this monomial, $d=0$ otherwise), and repeating this process with the new polynomial $r$ until we obtain either zero or a polynomial whose initial monomial is not divisible by any $in_\prec(g_j)$. Call the resulting polynomial a {\bf remainder} of $f$ with respect to $G$ and $\prec$.
This process is the same as the reduction algorithm found in \cite[Section 4]{Stokes}. It also agrees (up to the sign $(-1)^d$) with the usual reduction algorithm for elements of $\mathbb{F}[x_1,\ldots , x_n]$ found in, for example, \cite{Cox.Little.O'Shea}.

Define a {\bf left Gr\"obner basis} for a left ideal $I\subseteq \bigwedge$ to be a set $G=\{ g_1, \ldots , g_r \}\subset I$ such that $G$ left-generates $I$ and $\{in_\prec(g_1),\ldots , in_\prec(g_r)\}$ left-generates $in_\prec(I)$. This condition is straightforwardly equivalent to requiring that $f\in \bigwedge$ reduces to $0$ with respect to $G$ and $\prec$ whenever $f\in I$. By \cite[Theorem 6.5]{Stokes}, this is equivalent to requiring that for any $g_i,g_j\in G$, all $T_{g_i, \mathbf{x}^{\mathbf{c}}}$ and all $\mathbf{x}^{\mathbf{e}} S_{g_i,g_j}$ ($\mathbf{x}^{\mathbf{e}}\in M_{\bigwedge}$) reduce to zero with respect to $G$ and $\prec$. Using this condition \cite[Theorem 6.6]{Stokes} gives an algorithm, similar to Buchberger's algorithm, for extending a given generating set of $I$ to a left Gr\"obner basis of $I$ (\cite{Stokes} calls this a \emph{Gr\"obner left ideal basis}). The ideals we consider are homogeneous, so they are two-sided and we may drop the term ``left''.

\noindent\emph{Proof of Theorem~\ref{theorem:main}.} If $\prec_c$ is a coherent term order, it must satisfy $C_N$, so $B_i \prec_c A_i$ for some $1\le i\le N$. This implies $\mathbf{x}^{\mathbf{a}_i} \in in_{\prec_c}(I)$. It therefore suffices to show that no element of a generating set of $in_{\prec}(I)$ divides any $\mathbf{x}^{\mathbf{a}_j}, 1\le j\le N$. To this end, we use the algorithm of \cite[Theorem 6.6]{Stokes} to extend $H=\{\mathbf{x}^{\mathbf{b}_i}-\mathbf{x}^{\mathbf{a}_i}:1\le i\le N\}$ to a Gr\"obner basis $G=\{g_1,\ldots,g_r\}$ of $I$ with respect to $\prec$. We will show the elements of $G\setminus H$ involve monomials only of the form $\pm \mathbf{x}^{\mathbf{c}}\mathbf{x}^{\mathbf{a}_i}$ where $1\neq \mathbf{x}^{\mathbf{c}} \in M_{\bigwedge}$. This will suffice since $in_\prec(I)=\langle in_\prec(g_j) : g_j\in G\rangle$ and by the antichain condition, none of $\mathbf{x}^{\mathbf{b}_i}$, $\mathbf{x}^{\mathbf{c}}\mathbf{x}^{\mathbf{a}_i}$ (where $\mathbf{x}^{\mathbf{c}}\neq 1$) divide any $\mathbf{x}^{\mathbf{a}_j}, 1\le j\le N$.

Let $g\in H$ and $\mathbf{x}^{\mathbf{c}}\in M_{\bigwedge}$ with $\mathbf{x}^{\mathbf{c}}in_\prec(g)=0$. Then $\mathbf{x}^{\mathbf{c}}\neq 1$ and $T_{g,\mathbf{x}^{\mathbf{c}}} = \mathbf{x}^{\mathbf{c}}\mathbf{x}^{\mathbf{a}_i}$ (possibly equal to zero). So the remainder on reducing $T_{g,\mathbf{x}^{\mathbf{c}}}$ with respect to $H$ and remainders of other $T$-polynomials is either zero or also a monomial of this form. Any $T$-polynomial of a remainder of a $T$-polynomial is zero. 

Let $H'$ be the union of $H$ and the set of nonzero remainders of $T$-polynomials. 
The $S$-polynomial of two elements of $H'\setminus H$ is zero. Let $g_1\in H, g_2\in H'\setminus H, \mathbf{x}^{\mathbf{e}}\in M_{\bigwedge}$. By the antichain condition $g_2$ does not divide $in_{\prec}(g_1)$, so $\mathbf{x}^{\mathbf{e}}S_{g_1,g_2}$ is either zero or a monomial of the form $\pm \mathbf{x}^{\mathbf{c}}\mathbf{x}^{\mathbf{a}_i}$ for some $1\neq \mathbf{x}^{\mathbf{c}}\in M_{\bigwedge}$. This reduces similarly to the $T$-polynomials. 

Let $g_1,g_2\in H, \mathbf{x}^{\mathbf{e}}\in M_{\bigwedge}$. Then $S_{g_1,g_2}$ has the form $(-1)^{d_1}\mathbf{x}^{\mathbf{c}_1}\mathbf{x}^{\mathbf{a}_j}-(-1)^{d_2}\mathbf{x}^{\mathbf{c}_2}\mathbf{x}^{\mathbf{a}_i}$ 
(one or both of these terms may be zero), where $\mathbf{x}^{\mathbf{c}_1},\mathbf{x}^{\mathbf{c}_2}\in M_{\bigwedge}$. By the antichain condition $\mathbf{x}^{\mathbf{b}_i} \nmid \mathbf{x}^{\mathbf{b}_j}$ and $\mathbf{x}^{\mathbf{b}_j} \nmid \mathbf{x}^{\mathbf{b}_i}$, so $\mathbf{x}^{\mathbf{c}_1} \neq 1$ and $\mathbf{x}^{\mathbf{c}_2} \neq 1$. Thus reducing $\mathbf{x}^{\mathbf{e}}S_{g_1,g_2}$ with respect to $H'$ and nonzero remainders of other $S$-polynomials yields either $0$, or a monomial or binomial with term(s) of the form $\pm \mathbf{x}^{\mathbf{c}}\mathbf{x}^{\mathbf{a}_i}$, where $1\neq \mathbf{x}^{\mathbf{c}}\in M_{\bigwedge}$. The $S$-polynomial of any two such nonzero remainders, or such a nonzero remainder and an element of $H'$, is either zero or a monomial or binomial with terms of the form $\pm \mathbf{x}^{\mathbf{c}}\mathbf{x}^{\mathbf{a}_i}$ for some $1\neq \mathbf{x}^{\mathbf{c}}\in M_{\bigwedge}$, and the product of $\mathbf{x}^{\mathbf{e}}$ and this $S$-polynomial reduces similarly. \qed

\end{document}